\documentclass[11pt]{amsart}
\usepackage{amsbsy}
\begin{document}

%%%%%%%%%%%%%%%%%%%%%%%%%%%%%%%%%%%%%%%%%%%%%%%%%%%%%%%%%%%%%%%%%%%%
% Theorem, definition, lemma, proposition, corollary and proof
%%%%%%%%%%%%%%%%%%%%%%%%%%%%%%%%%%%%%%%%%%%%%%%%%%%%%%%%%%%%%%%%%%%%
%%%%%%%%%%%%%%%%%%%%%%%%%%%%%%%%%%%%%%%%%%%%%%%%
\newtheorem{theorem}{Theorem}%[section]
\newtheorem{proposition}{Proposition}%[section]
\newtheorem{lemma}{Lemma}%[section]
\newtheorem{corollary}{Corollary}%[section]%%
\newtheorem{definition}{Definition}%[section]
\newtheorem{remark}{Remark}%[section]
%%%%%%%%%%%%%%%%%%%%%%%%%%%%%%%%%%%
%%%%%%%%%%%%%%%%%%%%%%%%%%%%%%%%%%%%%%%%%%%%%%                  NEW
%%\newcommand{\be}{\begin{equation}}
%%\newcommand{\ee}{\end{equation}}
%%%%%%%%%%%%%%%%%%%%%%%%%%%%%%%%%%%%%%%%%%%%%%%
%%%%%%%%%%%%%%%%%%%%%%%%%%%%%%%%%%%%%%%%%%%%
\newcommand{\tex}{\textstyle}
%%\DeclareMathOperator{\Ind}{Ind}
%% \DeclareMathOperator{\Card}{Card}
%% \DeclareMathOperator{\Deg}{Deg}
%% \DeclareMathOperator{\dist}{dist}
%% \DeclareMathOperator{\Signature}{Signature}
%% \DeclareMathOperator{\MC}{MC}
%% \DeclareMathOperator{\sign}{sign}
%%  \DeclareMathOperator{\Int}{int}
%%%%%%%%%%%%%%%%%%%%%%%%%%%%%%%%%%%%%%%%%%%%
%%%%%%%%%%%%%%%%%%%%%%%%%%%%%%%%%%%%%%%%%%%%
\numberwithin{equation}{section} \numberwithin{theorem}{section}
\numberwithin{proposition}{section} \numberwithin{lemma}{section}
\numberwithin{corollary}{section}
\numberwithin{definition}{section} \numberwithin{remark}{section}
%%%%%%%%%%%%%%%%%%%%%%%%%%%%%%%%%%%%%%%%%%%%
%%%%%%%%%%%%%%%%%%%%%%%%%%%%%%%%%%%%%%%%%%%%
\newcommand{\ren}{\mathbb{R}^N}
\newcommand{\re}{\mathbb{R}}
\newcommand{\n}{\nabla}
\newcommand{\p}{\partial}
\newcommand{\iy}{\infty}
\newcommand{\pa}{\partial}
\newcommand{\fp}{\noindent}
\newcommand{\ms}{\medskip\vskip-.1cm}
\newcommand{\mpb}{\medskip}
%%%%%%%%%%%%%%%%%%%%%%%%%%%%%%%%%%%%%%%%%%%%%%%%%
\newcommand{\AAA}{{\bf A}}
\newcommand{\BB}{{\bf B}}
\newcommand{\CC}{{\bf C}}
\newcommand{\DD}{{\bf D}}
\newcommand{\EE}{{\bf E}}
\newcommand{\FF}{{\bf F}}
\newcommand{\GG}{{\bf G}}
\newcommand{\oo}{{\mathbf \omega}}
\newcommand{\Am}{{\bf A}_{2m}}
\newcommand{\CCC}{{\mathbf  C}}
\newcommand{\II}{{\mathrm{Im}}\,}
\newcommand{\RR}{{\mathrm{Re}}\,}
\newcommand{\eee}{{\mathrm  e}}
%%%%%%%%%%%%%%%%%%%%%%%%%%%%%%%%%%%%%%%%%%%%%%%%%%%%%%%%%%%%%%%%%%%%%%% L^2\rho...
\newcommand{\LL}{L^2_\rho(\ren)}
\newcommand{\LLL}{L^2_{\rho^*}(\ren)}
%%%%%%%%%%%%%%%%%%%%%%%%%%%%%%%%%%
%%%%%%%%%%%%%%%%%%%%%%%%%%%%%%%%%%%%%%%%%%%%%%%%%%%%
\renewcommand{\a}{\alpha}
\renewcommand{\b}{\beta}
\newcommand{\g}{\gamma}
\newcommand{\G}{\Gamma}
\renewcommand{\d}{\delta}
\newcommand{\D}{\Delta}
\newcommand{\e}{\varepsilon}
\newcommand{\var}{\varphi}
\newcommand{\lll}{\l}
\renewcommand{\l}{\lambda}
\renewcommand{\o}{\omega}
\renewcommand{\O}{\Omega}
\newcommand{\s}{\sigma}
\renewcommand{\t}{\tau}
\renewcommand{\th}{\theta}
\newcommand{\z}{\zeta}
\newcommand{\wx}{\widetilde x}
\newcommand{\wt}{\widetilde t}
\newcommand{\noi}{\noindent}
 %%%%%%%%%%%%%%%%%%%%%%%%%%%%%%%%%%%%%%%%%%%
\newcommand{\uu}{{\bf u}}
\newcommand{\xx}{{\bf x}}
\newcommand{\yy}{{\bf y}}
\newcommand{\zz}{{\bf z}}
\newcommand{\aaa}{{\bf a}}
\newcommand{\cc}{{\bf c}}
\newcommand{\jj}{{\bf j}}
\newcommand{\ggg}{{\bf g}}
\newcommand{\UU}{{\bf U}}
\newcommand{\YY}{{\bf Y}}
\newcommand{\HH}{{\bf H}}
\newcommand{\GGG}{{\bf G}}
\newcommand{\VV}{{\bf V}}
\newcommand{\ww}{{\bf w}}
\newcommand{\vv}{{\bf v}}
\newcommand{\hh}{{\bf h}}
\newcommand{\di}{{\rm div}\,}
\newcommand{\ii}{{\rm i}\,}
%%%%%%%%%%%%%%%%%%%%%%%%%%%%%%%%%%
%%%%%%%%%%%%%%%%%%%%%%%%%%%%%%%%%%%%%   VAG, NEW
\newcommand{\inA}{\quad \mbox{in} \quad \ren \times \re_+}
\newcommand{\inB}{\quad \mbox{in} \quad}
\newcommand{\inC}{\quad \mbox{in} \quad \re \times \re_+}
\newcommand{\inD}{\quad \mbox{in} \quad \re}
\newcommand{\forA}{\quad \mbox{for} \quad}
\newcommand{\whereA}{,\quad \mbox{where} \quad}
\newcommand{\asA}{\quad \mbox{as} \quad}
\newcommand{\andA}{\quad \mbox{and} \quad}
\newcommand{\withA}{,\quad \mbox{with} \quad}
\newcommand{\orA}{,\quad \mbox{or} \quad}
\newcommand{\atA}{\quad \mbox{at} \quad}
\newcommand{\onA}{\quad \mbox{on} \quad}
\newcommand{\ef}{\eqref}
\newcommand{\mc}{\mathcal}
\newcommand{\mf}{\mathfrak}

\newcommand{\ssk}{\smallskip}
\newcommand{\LongA}{\quad \Longrightarrow \quad}
%%%%%%%%%%%%%%%%%%%%%%%%%%%%%%%%
%%%%%%%%%%%%%%%%%%%%%%%%%%%%%%%%%%
\def\com#1{\fbox{\parbox{6in}{\texttt{#1}}}}
%%%%%%%%%%%%%%%%%%%%%%%%%%%%%%%%%%
%%%%%%%%%%%%%%%%%%% From Paper1
\def\N{{\mathbb N}}
\def\A{{\cal A}}
\newcommand{\de}{\,d}
\newcommand{\eps}{\varepsilon}
\newcommand{\be}{\begin{equation}}
\newcommand{\ee}{\end{equation}}
\newcommand{\spt}{{\mbox spt}}
\newcommand{\ind}{{\mbox ind}}
\newcommand{\supp}{{\mbox supp}}
\newcommand{\dip}{\displaystyle}
\newcommand{\prt}{\partial}
\renewcommand{\theequation}{\thesection.\arabic{equation}}
\renewcommand{\baselinestretch}{1.1}
%%%%%%%%%%%%%%%%%%%%%%%%%%%%%%%%%%%%%%%%%%%%%%%
\newcommand{\Dm}{(-\D)^m}

%%%%%%%%%%%%%%%%%%%%%%%%% VICTOR
\title
%%%%%[Positivity] %%%%%%%%%%%%%%%%%%%%%%%%%
 {\bf %% An analytic regularization
 Well-posedness of the Cauchy
problem for a fourth-order thin film equation via regularization approaches}

\author{Pablo~\'Alvarez-Caudevilla\hspace{2mm}
\vspace{1mm}\\
{\small Departamento de Matem{\'a}ticas},\\ {\small U.~Carlos III de Madrid},
{\small 28911 Legan{\'e}s, Spain}\\
{\small e-mail: pacaudev@math.uc3m.es}\vspace{1mm}\\
Victor~A.~Galaktionov \hspace{2mm}
\vspace{1mm}\\
{\small Department of Mathematical Sciences},\\ {\small University of Bath},
{\small Bath BA2 7AY, UK}\\
{\small e-mail: vag@maths.bath.ac.uk}\vspace{1mm}}


\thanks{This work has been partially supported by the Ministry of Economy and Competitiveness of
Spain under research project MTM2012-33258.}

\date{\today}

%%%%%%%%%%%%%%%%%%%%%%%%%%%%%%%%%%%%%%%%%%%%%%%%%%%%%%%%

%%\setlength{\topmargin}{5mm} \pagestyle{myheadings} \markboth{}
%%{Fourth-order thin film equation}

%%\begin{document}

%%%%$$ $$

\begin{abstract}

This paper is devoted to some aspects of well-posedness
 of the
Cauchy problem (the CP, for short) for a quasilinear degenerate
fourth-order parabolic \emph{thin film equation} (the TFE--4)
 \be
 \label{0.1}
 u_{t} = -\nabla \cdot(|u|^{n} \nabla
\D u) \inB \ren \times \re_+, \quad u(x,0)=u_0(x) \inB \ren,
 \ee
 where $n>0$ is  a fixed exponent,
 with bounded smooth compactly supported initial data. Dealing with the CP (for, at least,
 $n \in (0, \frac 32)$) requires
 introducing
classes of infinitely changing sign solutions that are oscillatory
close to finite interfaces.
 The main goal of the paper is to
detect proper solutions of the CP for the degenerate TFE--4 by
uniformly parabolic analytic $\e$-regularizations at least for
values of the parameter $n$ sufficiently close to 0.

Firstly,  we study an analytic ``homotopy" approach based on {\em
a priori} estimates for solutions of uniformly parabolic analytic
$\e$-regularization problems of the form
  $$
   u_{t} = -\nabla \cdot
(\phi_\e(u)\nabla \D u) \inB \ren \times
    \re_+,$$
where $\phi_\e(u)$ for $\e \in (0,1]$ is an analytic
$\e$-regularization of the problem \eqref{0.1}, such that
$\phi_0(u)=|u|^n$ and $\phi_1(u)=1$, using a more standard classic technique of
passing to the limit in integral identities for weak solutions.
However, this argument has been demonstrated  to be non-conclusive;
basically due to the lack of a complete optimal estimate-regularity
theory for these types of problems.

Secondly, to resolve that issue more successfully, we study a more
general similar
  analytic ``homotopy
transformation" in both the parameters, as $\e \to 0^+$ and $n \to
0^+$, and
 describe
 {\em branching} of solutions of the TFE--4 from the solutions of the notorious \emph{bi-harmonic
equation}
  $$
  u_t=-\D^2 u  \inB \ren \times \re, \quad u(x,0)=u_0(x) \inB
  \ren,
$$
 which
 describes
some qualitative oscillatory properties of CP-solutions of
 \ef{0.1} for small $n>0$ providing us with the uniqueness of solutions for the problem \eqref{0.1}
 when $n$ is  close to 0.

 Finally,  {\em
 Riemann-like problems} occurring in a boundary layer
 construction, that occur close to nodal sets of the solutions, as $\e \to 0^+$,
    are discussed in other to get uniqueness results for the TFE--4 \eqref{0.1}. 

\end{abstract}

%%%%%%%%%%%%%%%%%%%%%%%%%%%
\maketitle

%%%%%%%%%%%%%%%%%%%%%%%%%%%%%%%%%%%%%%%%%%%%%%%%
\section{Introduction:  the Cauchy and free boundary problems for the TFE--4}
 \label{S1}

\subsection{Main model and their applications}

\noindent In this paper, we study some aspects  of well-posedness
of the Cauchy problem (the CP) for a nowadays well-known
fourth-order quasilinear evolution equation of parabolic type,
called the {\em thin film equation} (the TFE--4), with an exponent
$n>0$,
\begin{equation}
\label{i1}
    u_{t} = -\nabla \cdot(|u|^{n} \nabla \D u)
 \quad \mbox{in} \quad \ren \times \re_+\,,
 \quad u(x,0)=u_0(x) \inB \re^N,
\end{equation}
with  bounded, sufficiently smooth, and compactly supported initial
data $u_0$ (not necessarily positive)
with an arbitrary dimension  $N \ge 1$. Note that these initial conditions could be relaxed, (for example $u_0\in L^1\cap L^\infty$) 
however, it is not the purpose of this 
work to analyse the problem from the perspective of different possible initial conditions.

\ssk

Equation \eqref{i1} arises in numerous physical related areas with
applications in thin film, lubrication theory, and in several
other hydrodynamic-type problems. In particular, those equations
model the dynamics of a thin film of viscous fluid, as the
spreading of a liquid film along a surface, where $u$ stands for
the height of the film. Then clearly assuming $u \ge 0$ naturally
leads to a \emph{free boundary problem} (an FBP) setting; see
below. Specifically, when $n=3$, we are dealing with a problem in
the context of lubrication theory for thin viscous films that are
driven by surface tension and when $n=1$ with Hele--Shaw flows.
However, in this work, we are considering solutions of changing
sign. Such solutions can have some biological motivations
\cite{KingPers}, to say nothing of general PDE theory, where the
CP-settings were always key.

%%%%%%%%%%%%%%%%%%%%%%%%%%%%%%%%%%%%%%%%%%%%%%%%%%%%%%%%

\subsection{Main results}

%%%%%%%%%%%%%%%%%%%%%%%%%%%%%%%%%%%%%%%%%%%%%%%%%%%%%%

As a more successful approach, among others, to clarify the
well-posedness of the CP, we perform an analytic ``homotopic"
approach 
%when the parameter $n$ goes to zero 
from our original
equation \eqref{i1} to an equation from which we can extract
information about the solutions of equation \eqref{i1},
% at least when the parameter $n$ is sufficiently close to zero. 

Namely, we
develop a {\em homotopic
deformation} 
 from the TFE--4 \eqref{i1} to
the classic and well-known {\em bi-harmonic equation}
\begin{equation}
\label{s1}
 \tex{
    u_{t} = -\D^2 u\quad \hbox{in} \quad \ren \times \re_+\,,
    \quad u(x,0)=u_0(x) \inB \ren.
    }
\end{equation}
It is well-known that,
for any smooth compactly supported data $u_0$, satisfying the
natural ``growth condition at infinity"
\begin{equation*}
    u_0 \in L_{\rho^*}^2(\ren), \quad \mbox{where} \quad \rho^*(y)={\mathrm e}^{-a |y|^{4/3}},
  \quad a={\rm const.}>0 \,\,\,\mbox{small},
\end{equation*}
the {\em bi-harmonic} equation \eqref{s1} admits a unique classic solution given by the convolution
Poisson-type integral,
\begin{equation}
\label{s2}
 \tex{
    \tilde{u}(x,t)=b(x,t)\, * \, u_0(x) \equiv t^{-\frac N4} \int\limits_{\ren} F((x-z)t^{-\frac 14})
     u_0(z)\, {\mathrm d}z,
    }
\end{equation}
 where $b(x,t)$ is the fundamental solution
\begin{equation}
\label{s3}
 \tex{
    b(x,t)=t^{-\frac N 4} F(y), \quad y=\frac{x}{t^{1/4}},
    }
\end{equation}
of the operator $\frac{\p}{\p t} + \D^2$. The oscillatory rescaled
kernel $F(y)$ is the unique solution of the linear elliptic
problem
\begin{equation}
\label{s4}
 \tex{
    {\bf B}F \equiv -\D_y^2 F + \frac{1}{4}\, y \cdot \nabla_y F  +\frac{N}{4}\, F=0
    \quad \hbox{in} \quad \ren,\quad \int_{\ren} F(y) \, {\mathrm
    d}y=1.
    }
\end{equation}

\ssk

\noi{\sc \underline{First regularization.}} Thus, firstly, 
%to this end, 
we perform a 
%simpler 
homotopic deformation
assuming that $n>0$ is a sufficiently small and {\em fixed}
exponent. We are then actually talking  about
 some  ``homotopic classes" (understood here not in the classic sense from degree operator theory)
  of degenerate parabolic PDEs. 
  
  More precisely, we say that
 the TFE \eqref{i1} is ``homotopic to the linear PDE \eqref{s1}"
if there exists a family of uniformly parabolic equations (a {\em
homotopic deformation}) with a coefficient 
$$\phi_\e(u)>0\quad \hbox{analytic in both variables $u\in\re$ and $\e\in (0, 1]$},$$ 
with unique
analytic solutions $u_\e(x,t)$ of the problem
\begin{equation}
\label{cp4}
    u_{t} = -\nabla \cdot (\phi_\e(u)\nabla \D u) \inB \ren \times
    \re_+, \quad u(x,0)=u_0(x) \inB \ren, 
    %%\,\,\, \mbox{such that}
\end{equation}
such that
\begin{equation}
\label{cp5}
 \phi_1(u)=1, \andA
    \phi_\e(u)\rightarrow |u|^n\quad \hbox{as}\quad \e \rightarrow 0^+ \quad
    \hbox{uniformly on compact subsets}\,,
\end{equation}
so that $u(x,t)$ can be approximated by $u_\e(x,t)$ as $\e \to
0^+$.

A possible and quite natural
homotopic path (to be used in this work) is
\begin{equation}
 \label{cp51}
    \phi_\e(u):=\e^n+(1-\e)(\e^2+u^2)^{\frac n2}, \quad \e \in (0, 1].
\end{equation}
Then, indeed,  the non-degenerate uniformly parabolic equation
\eqref{cp4} admits a unique (at least, locally in time) classic
solution $u=u_\e(x,t)$, which is an analytic function in all the
three variables $x$, $t$, and $\e$, just using classic parabolic
theory \cite{EidSys, Fr}, for any $\e \in (0,1]$.

Therefore, the
homotopic deformation described above is basically a continuous
deformation from the TFE--4 to the bi-harmonic equation through
the non-degenerate equation \eqref{cp4}, \eqref{cp5} for $\e\in
(0,1]$, for which we also know key features, such as the uniqueness of solutions, 
oscillatory changing sign properties of them and the well-posedness of the Cauchy problem.

%%\ssk

Hence, given the well-defined  analytic functional
family  (a curve or a path),
 \be
 \label{Fem}
 {\mathcal P_\phi}=
\{u_\e(x,t)\}_{\e \in (0,1]},
 \ee
 as usual in extended semigroup theory,
  a solution $u(x,t)$ of the CP for the TFE--4 \ef{i1}
is then a function, satisfying in a
standard pointwise (or even uniformly, if better estimates are
available), or in a weaker sense,
 \be
 \label{Fem1}
 u_\e(x,t) \to u(x,t) \asA \e \to 0^+.
  \ee
Consequently, as a first step we perform such a homotopic deformation from the TFE--4 
\eqref{i1} to the {\em bi-harmonic equation} \eqref{s1}, through
the non-degenerate $\e$-regularization equation \eqref{cp4}, \eqref{cp5}, using energy methods as those
used by Bernis--Friedman in their pioneering work \cite{BF1}.

  However, as we will see, passing to the limit in some more or
  less standard integral identities for $\{u_\e\}$ is not
  sufficient to convincingly distinguish the CP and FBP solutions
  (both can admit the same  weak setting). 
  
  In other words, the weak
formulation of the CP does not allow one to clarify key features
of the CP-solutions. So, in this process, open questions remain,
and this standard analysis does not solve the problem of
uniqueness of solutions of the CP.

  Indeed, in general, under necessary natural estimates on $\{u_\e\}$, deformation \ef{Fem1}
  is able to detect {\em a solution} of the CP, and does not
  guarantee its {\em conventional uniqueness} via existence of a limit in \ef{Fem1}.
   So that, \ef{Fem1} may contain several (or even infinitely many)
  {\em partial limits}, since the end point $\e=0^+$ of the path \ef{Fem} is singular; see Section\;\ref{S2}. 

  Moreover, nothing will also guarantee that a proper CP-solution
   {\em is actually independent of the type
  of an analytic $\e$-regularization $\phi_\e(u)$ in the path $\ef{Fem}$.} This is the strongest type of uniqueness,
   which is called the {\em absolute uniqueness}. In fact,
   this is very difficult
  to achieve for higher-order PDEs\footnote{It seems that, in a most general setting, such absolute uniqueness
  results are non-existent in principle (though these could be useless if they exist).}.

 Note that, even for second-order parabolic equations, such a type of  uniqueness of {\em proper
  solutions} (extremal: e.g., {\em minimal} ones in blow-up problems), which is
  guaranteed by the Maximum Principle and comparison techniques, was proved
   under the assumption on the {\em monotonicity in $\e$ of any
  $\e$-regularization applied} (then any of such regularization leads to the unique proper minimal solution); see \cite[Ch.~7]{GalGeom}.
  Hence, for the TFE--4, we do not stress  upon such a
  strong uniqueness issue.

\ssk

\noi{\sc \underline{Second regularization: a double limit.}}
  Therefore, secondly, our further improved
  definition of the proper CP-solutions assumes also a second
  limit as $n \to 0^+$ (with a simpler regularization; see below), i.e., a continuous connection with
  solutions of the {\em bi-harmonic equation} \eqref{s1} which we know explicitly. Thus, $u(x,t)$ in \ef{Fem1} is a CP-solution
  of the TFE--4 \ef{i1} if
   \be
   \label{Fem11}
   u(x,t) \to \tilde u(x,t) \asA n \to 0^+,
    \ee
where $\tilde u(x,t)$ is denoted by \eqref{s2}.

    Actually, instead of independent limits \ef{Fem1}, \ef{Fem11},  we will need to perform
a kind of {\em double} limit as  $\e,\, n \to 0^+$, where special
restrictions on two parameters will be required. 

However, as discussed at the end of Section\;\ref{S2} even by
this double limit  as  $\e,\, n \to 0^+$, using standard integral
identities, it is not still possible to identify the problem which
the limit belongs to. Hence, something different must be done.

Additionally open questions remain  regarding regularity of the
solutions for the TFE--4 \eqref{i1}, especially when $N\geq 2$.
Note that, for $N=1$, a complete analysis with several important
uniform estimates  was already done in \cite{BF1}. This issue is
essentially the reason why the energy methods do not work in
$\re^N$ with $N\geq 2$ since extensions of that analysis are not straightforward. 

Consequently, as a different and simpler, but more general,
regularizing version, a homotopic path as $n \to 0^+$ and $\e \to
0^+$
  can be used separately for a
 proper definition of the solutions of the Cauchy problem for the
 TFE--4 (\ref{i1}). However, imposing certain relation between the parameters $\e$ and $n$.

 Hence, using
the following simpler approximation in \ef{cp4}:
  \be
  \label{phi2}
  \phi_\e(u)=(\e^2+u^2)^\frac n2,
   \ee
 this assumes that a correct solution of the CP for the TFE--4 \ef{i1} is  that one,
  which can be continuously
   deformed via
   this ``double" limit to the unique solution of the {\em bi-harmonic equation} \eqref{s1}.

  Now, we can state  the main result of this paper, which provides us with a definition of the proper solution of the Cauchy problem \eqref{i1}
   as well as the uniqueness of solutions and well-posedness of the CP, when the parameter $n$ is sufficiently close to zero.

   \begin{theorem}
\label{principal} Assume the condition
 \be
 \label{lim12}
n \, |\ln \e(n)| \to 0 \asA n \to 0^+,
 \ee
 holds and the regularization
family $\{u_\e\}$ is uniformly bounded. Then:

{\rm (i)}  The solution of the perturbed problem
\begin{equation*}
  u_t=- \n \cdot ((\e^2+u^2)^{\frac n2} \n \D u), \quad
u(x,0)=u_0(x),
 \end{equation*}
converges uniformly to the  solution \ef{s2} of the Cauchy problem
of the bi-harmonic equation
\begin{equation*}
  u_t=- \D^2 u, \quad
u(x,0)=u_0(x),
 \end{equation*}
in
 the  limit $\e=\e(n) \to 0^+$ and $n \to 0^+$, which is unique.

{\rm (ii)}  Provided that the convolution
$$\tex{\varphi_1(t)= -  \int_0^t  \n b(t-s) *  \ln | \tilde u(s)|\n \D \tilde u(s) \, {\mathrm d} s,}$$ remains
uniformly bounded for the solution \ef{s2}, the rate of
convergence as $n \to 0$ of the $($formal$)$ asymptotic expansion $u = \tilde{u} + V$ is given by
$$ \tex{\quad V:= n \varphi_1 + o(n),\quad \hbox{with}\quad \varphi_1(t)= -  \int_0^t  \n b(t-s) *  \ln | \tilde u(s)|\n \D \tilde u(s) \, {\mathrm d} s.}$$
\end{theorem}

Of course, a similar result can be more easily established for the
previous, single-limit path \ef{cp51}, but the present
double-limit one is indeed more general, though assumes the
necessary restriction \ef{lim12} on parameters. Remember that this restriction is crucial to have such results.

Consequently, thanks to the previous result, we can assure that there
exists a branch of solutions of the TFE-4 \eqref{i1} emanating at
$n=0^+$ from the unique solution of the parabolic bi-harmonic
equation \eqref{s1}, deforming the solutions of the TFE-4
\eqref{i1} as $n \to 0^+$ and $\e\to 0^+$ via an analytic path (through the $\e$ regularization \eqref{phi2}) and inheriting the
oscillatory and changing sign properties of the linear flow.

Also, since the parabolic bi-harmonic equation \eqref{s1} is well-posed, through the homotopy 
deformation performed by Theorem\;\ref{principal} we can achieve 
such a well-posedness for the CP \eqref{i1}, at least when the parameter $n$ is sufficiently close to zero.

Note that we are going to develop a
branching theory for difficult {\em nonlinear degenerate parabolic
PDEs}, requiring special relationships between $n \to 0^+$
and the regularization parameter $\e=\e(n) \to 0$. Without such a
relation, no proper limits can be traced out at all.

In the final section, a further and deeper study of these limits (particularly, for
establishing a proper uniqueness of CP-solutions of the TFE--4 via
regularization) leads to difficult
 {\em boundary layer}-type problems that remain open in a
 sufficient generality.
 To be specific,  we concentrate here on the analysis of the problematic limit as $\e\to
 0^+$ for \ef{phi2}
  in order to prove the uniqueness of solutions
 for the TFE--4 \eqref{i1}. 
 
 Although, we perform our analysis in arbitrary dimensions, the general (say, the absolute one)  uniqueness
 result is still an
 open problem and probably non-achievable in principle.
 To do so, we need to understand the zero structure sets of the solutions, especially close to the interfaces. This leads us to analyze the {\em boundary layers}
 that occur close to those nodal sets. In other words, to study an asymptotic problem near the interfaces, usually called {\em Riemann's problems} for the TFE--4.
Indeed, the well-posedness of the limiting problem when $\e \to 0^+$ will strongly depend on the asymptotics of the solution for the TFE--4 \eqref{i1}
close to the interfaces, as discussed in the final section of this paper.

%%%%%%%%%%%%%%%%%%%%%%%%%%%%%%%%%

\subsection{FBP and CP settings: similarities and distinctive features}

%%%%%%%%%%%%%%%%%%%%%%%%%%%%%%%%%%%%%%

 First serious attempts to create a proper mathematical theory for
 TFEs--4  were made in the 1980s. Since then, hundreds of
 papers and a number of monographs have dealt seriously with such
 higher-order parabolic PDEs.
 We refer to e.g. \cite{EGK1, EGK2, Gia08, Grun04} for
most recent  surveys and for extended lists of references
concerning physical derivations of various models, key
mathematical results, and further applications.  Nowadays (and at
least  since the 1980s), such equations  play quite a  special and,
even a key role, in general nonlinear PDE theory.

For this particular work  it is crucial to note that the TFE--4 \eqref{i1} is written for solutions of
changing sign, which, as we show,  can occur in the CP and also in
some free-boundary problems.

As mentioned above, since the 1980s,
it has been customary to consider {\em non-negative solutions} of
the TFE--4, and moreover
 the existing mathematical PDE theory was created mostly for
 solutions 
 $$u=u(x,t) \ge 0.$$ 
  
  However,  {\em solutions of changing
sign} have been under scrutiny already for a few years;
see \cite{BW06, EGK2, EGK4} and references therein.
To deal with such a ``dichotomy":
   $$\hbox{ {\em non-negative} FBP-
or {\em infinitely sign-changing} CP-solutions},$$
  let us recall
our main convention (to be discussed)
 \be
 \label{Dich}
 \begin{matrix}
 \mbox{oscillatory sign changing solutions are related to the CP,
 while} \ssk \ssk \\
 \mbox{non-negative solutions appear for the standard FBP (as in \cite{BF1,BHQ})},
 \end{matrix}
 \ee
 for not that large  $n>0$, actually,
$$\hbox{for any}\quad 0<n < n_{\rm h}= 1.7587... \, .
   $$
    See
\cite{EGK2} for any further details, where that upper bound is
obtained using numerical analysis, and \cite{PetI} for proper
estimates of this {\em heteroclinic critical exponent} $n_{\rm
h}$.

For both problems, the CP and FBP,
we assume that the solutions satisfy the following standard free boundary conditions:
\begin{equation}
\label{i3}
    \left\{\begin{array}{ll} u=0, & \hbox{zero-height,}\\
    \nabla u=0, & \hbox{zero contact angle,}\\
    -{\bf n} \cdot \nabla (|u|^{n}  \D u)=0, &
    \hbox{zero-flux (conservation of mass)}\end{array} \right.
\end{equation}
at the singularity surface (interface) $\Gamma_0[u]$, which is the
lateral boundary of the  support
\begin{equation*}
    \hbox{supp} \;u \subset \ren \times \re_+,
\end{equation*}
where ${\bf n}$ stands for the unit outward normal to
$\Gamma_0[u]$, which is assumed to be sufficiently smooth, at
least, a.e.\footnote{A very difficult result to prove, still a
fully open problem for any $N \ge 2$. It seems that possible types
of ``singular cusps" of interfaces occurring via their
self-focusing (e.g., the simplest Graveleau--Aronson-type radial
``filling a hole" in PME theory known for \ef{PME1} \cite{AA})
were not studied at all for \ef{i1}.} The condition of zero flux
might be interpreted as
\begin{equation*}
    \lim_{\hbox{dist}(x,\Gamma_0[u])\downarrow 0}
    -{\bf n} \cdot \nabla (|u|^{n}  \D u)=0,
\end{equation*}
 again for sufficiently smooth interfaces.
 Note that there can be other and different free-boundary conditions (e.g., without the zero contact angle
 one),
     and a lot of work has
 been directed to the discussion of that matter, for instance \cite{Otto}.

The problem is completed with bounded, smooth, integrable,
compactly supported initial data:
\begin{equation}
\label{i4}
    u(x,0)=u_0(x) \quad \hbox{in} \quad \Gamma_0[u] \cap \{t=0\}.
\end{equation}
Thus, regarding the initial data, unlike the FBP, in the CP for
\eqref{i1} in $\ren \times \re_+$, one needs to pose bounded
compactly supported initial data \eqref{i4} prescribed in the
whole $\ren$. 

Intuitively (and rather loosely speaking),  the
principal difference between the CP and  the FBP is that, for the
CP, {\em no free boundary conditions} should be prescribed {\em a
priori}, while, for the FBP, these must be clearly given as done
in \ef{i3}. 

In this connection, it is worth mentioning that it is
known that solutions of both the CP and the FBP can be obtained by
proper $\e$-regularizations, i.e., {\em without} specifying the
free-boundary conditions. 

What is a principal fact is that those
$\e$-regularizations for the CP and the FBP must be {\em
different}, unless both problems coincide. The latter  can
actually happen for the TFE-4 for larger $n \ge n_{\rm h}$,  when
solutions of the CP lose their infinitely oscillatory properties,
or even for $n \ge 2$ only. In general, the question of $n$'s, for
which ``the CP $=$ the FBP", remains open.

Note also that
 there exist many
other FBPs with slightly different conditions on the interface
(cf. e.g., {\em Stefan} and {\em Florin} classic FBPs for the heat
equation, \cite[Ch.~8]{GalGeom}).  Of course, for the CP, it would
be very interesting to derive the actual and sharp enough
conditions on the free boundary (of course, \ef{i3} are valid, but
these have nothing to do with the optimal regularity for the CP),
but this is difficult even in 1D, to say nothing of any proof; see
\cite[\S\,6]{EGK2}.

Obviously, due to the zero flux condition in \ef{i3}, for the FBP,
the total mass of solutions is preserved in time. For the CP, as
we know, this must be also true, but not that straightforward.
 Here, it is natural to declare another additional convention (to be
 treated as well): at least for $n \in (0, \frac 32)$,
  \be
  \label{Conv1}
  \mbox{the solutions of the CP are {\em smoother} at the interface than
  for the FBP.}
   \ee
Actually, we have got some evidence about that convention in \cite{EGK2}.  Indeed, close to the
interface for $t>0$, after an expected {\em ``interface waiting
time"}, in the sense of the ``multiple zeros at the interfaces", it follows that
 \be
 \label{Conv2}
  \fbox{$
 \mbox{CP:} \quad u(x, t) \in C^{[\frac 3n]-1}_x,
 $}
 \ee
where $[\cdot]$ denotes the integer part\footnote{Note that inner,
``transversal" zeros of solutions can be less regular, but these
are not that interesting, unlike the key interface ones. However,
``transversal" zeros are accumulated near interfaces, so
\ef{Conv2} should not be understood literally, and it just shows
the smoothness of the {\em non-oscillatory envelope} of solutions
near interfaces (see \cite[\S~7]{EGK2}), or derivatives are
assumed to be calculated over neighbourhoods without touching
enough those ``transversal" zeros.}. This means that, for
 small $n>0$, the CP-solutions can be arbitrarily smooth in $x$ at the
 interfaces. This will be dealt with rather seriously later.

 Thus,  under the same
zero flux condition at finite interfaces (which should be
established separately), the mass is preserved in the CP as well.
 Let
\begin{equation*}
 \tex{
    M(t):=\int u(x,t) \, {\mathrm d}x
    }
\end{equation*}
be the mass of the solution, where integration is performed over
smooth support ($\ren$ is allowed for the CP only). Then,
differentiating $M(t)$ with respect to $t$ and applying the
divergence theorem (under natural regularity assumptions on
solutions and free boundary), we have that
\begin{equation*}
 \tex{
  J(t):=  \frac{{\mathrm d}M}{{\mathrm d}t}= -
  \int\limits_{\Gamma_0\cap\{t\}}{\bf n} \cdot \nabla
     (|u|^{n}  \D u )\, .
     }
\end{equation*}
Then, the mass is conserved if
   $ J(t) \equiv 0$, which is assured by the flux condition in
   \eqref{i3}.

%%%%%%%%%%%%%%%%%%%%%%%%%%%%%%%%%

 %%\ssk
%%%%%%%%%%%%%%%%%%%%%%%%%%%%%%%%%%%%%%%%%%%%%%%%%%%%%%%%%%%%%%
\subsection{Solutions of changing sign and comparison with porous medium flows}
As mentioned above, most of the existing results for the TFE--4
deal with {\em non-negative solutions} with compact support of
various FBPs, which are often more physically relevant. In this
context, we should point out that such approximations for {\em
non-negative} and non-changing sign solutions,  with various
non-analytic (and non-smooth) regularizations, for example, of the
form 
$$\phi_\e(u)=|u|^n +\e,$$ 
which is not only non-analytic for $n<2$, but even
$\not \in C^2$, have been widely used before in TFE--FBP theory as
a key foundation; cf. \cite{BBP}, \cite{BF1}, and \cite{BGK}.

Indeed, as pointed out in \cite{BHQ}, for non-negative solutions,
when $n$ approximates $0^+$, the limiting problem is always a free
boundary problem. Owing to the oscillatory
 behaviour of the solutions of the CP for the linear bi-harmonic equation \eqref{s1},
  the CP cannot be the limiting problem
in this  case.

In other words, in general, positivity of solutions of some FBPs
for the TFE--4 and other related equations is achieved via some uniformly parabolic, but
sufficiently ``singular" (and surely non-analytic) $\e$-regularization of the
PDE. In fact, this creates special kinds of ``obstacle FBPs",
where the obstacle appears to get positive solutions.
This was first proved by specific singular regularizations
introduced in the seminal TFE-paper \cite{BF1}, devoted mainly to
FBP settings.

On the other hand, it turned out that  the classes of the
so-called ``oscillatory solutions of changing sign" of (\ref{i1})
were rather difficult to tackle rigorously by standard and classical
methods. Specifically, due to the fact that any kind of a detailed
analysis for higher-order equations is much more difficult than
for their second-order counterparts, such as the notorious classic
\emph{porous medium equation} (PME--2)
  \be
  \label{PME1}
   u_{t} = \D (|u|^{n-1}u)
\inA,
   \ee
   in view of the lack of the Maximum Principle, comparison methods,
order-preserving semigroups, and potential properties of the
operators involved. Thus, practically all the existing methods for
monotone or variational operators are not applicable for \ef{i1}.

The PME--2 \ef{PME1} can be interpreted as a nonlinear degenerate
   extension of the classic \emph{heat equation} for $n=0$,
   \be
   \label{HE1}
   u_t=\D u \inA.
    \ee
 Note that passing to the limit $n \to 0^+$ in the PME--2 \ef{PME1} for non-negative solutions
  used to be a difficult
 mathematical problem in the 1970s-80s, which exhibited typical
 features (but clearly simpler than those ones in the TFE case)
 of a ``homotopy" transformation of PDEs.
 This study was initiated by Kalashnikov in 1978
 \cite{Kal78} for the one-dimensional case. Further detailed results in $\ren$ were obtained in
 \cite{Ben81}; see also \cite{Cock99}. More recent involved
 estimates were obtained in \cite{Pan08, Pan07} for the 1D PME--2
 establishing the rate of convergence of solutions
 as $n \to 0^\pm$, such as $O(n)$ as $n \to 0^-$ (i.e, from $n<0$, the fast diffusion
range, where solutions are smoother)
  in $L^1(\re)$ \cite{Pan08}, and $O(n^2)$ as $n \to 0^+$ in
 $L^2(\re \times (0,T))$ \cite{Pan07}.

 However, most of such convergence
 results are also obtained for {\em non-negative} solutions of the
 PME--2.
 For solutions of changing sign, even for this second-order  PME--2 \ef{PME1},
 there are some open problems; see \cite{GHCo} for
 references and further details.

Thus,  in the twenty-first century, higher-order TFEs such as
\ef{i1}, though looking
  like a natural  and not-that-involved counterpart/extension of the PME--2 \ef{PME1},
   corresponding mathematical
TFE theory gets essentially more complicated with several key
problems remaining  open still (such as regularity, uniqueness in general, etc).

\vspace{0.2cm}

The outline of the paper is: in Section \ref{S2} we perform an homotopic approach passing to the limit in integral identities for weak solutions. 
In Section \ref{S3}, applying a \emph{branching argument} we study the homotopy \emph{double limit} as $n,\e \to 0^+$. Finally,  in Section \ref{S5}
the limit as $\e \to 0^+$ is analysed in order to prove the uniqueness of solutions for the TFE--4 \eqref{i1}.

%%%%%%%%%%%%%%%%%%%%%%%%%%%%%%%%%%%%%%%%%%%%%%%%%%%%%%%%%%
\section{The $\e$-regularization problem: passage to the limit in a weak sense}
 \label{S2}

%%%%%%%%%%%%%%%%%%%%%%%%%%%%
\subsection{Preliminary estimates}

For any $\e\in (0, 1]$, we denote by $u_\e(x, t)$ the unique
solution of the CP for the regularized non-degenerate uniformly
parabolic equation \eqref{cp4}, \ef{cp51} with the same data
$u_0$. 

By classic parabolic theory \cite{EidSys, Fr}, this family
$$\{u_\e(x,t)\}\quad \hbox{is continuous (and analytic) in $\e\in (0, 1]$},$$ 
in
any natural functional topology, at least, on a time interval
$[0,T]$. Indeed, we also have that all the derivatives are
H\"{o}lder continuous in $\overline{\O} \times [0,T]$. However, as we show below, the main problem is the behaviour as
$\e \rightarrow 0^+$, where the regularized PDE loses its uniform
parabolicity.

Here, by
$\O$ we denote either $\ren$, or, equivalently, the bounded domain
${\rm supp} \,u \cap\{t\}$, i.e., the section of the support.

Note again that,  for second-order parabolic equations with the
Maximum Principle, such regularization-continuity approaches are
typical for constructing unique proper extremal  (say, minimal)
solutions with various singularities (e.g., finite time blow-up,
extinction, finite interfaces, etc.); see \cite[Ch.~4-7]{GalGeom}
as a source of key references and basic results.

However, for
higher-order degenerate parabolic flows admitting strongly
oscillatory solutions of changing sign, such a homotopy-continuity
approach generates a number of difficult problems. In fact,
despite the fact that the passage to the limit as $\e \to 0^+$ looks
like a reasonable way to define a proper solution of the TFE, we
expect that there are always special classes of compactly
supported initial data, for which such a limit is non-existent
and, moreover, there are many partial limits, thus defining a
variety of different solutions (meaning non-uniqueness).

Thus, to study such a limit \ef{Fem1} for the problem \eqref{cp4},
\ef{cp51}, when $\e \rightarrow 0^+$, we firstly need some
estimates for its regularized solutions $\{u_\e(x, t)\}$.

\begin{proposition}
\label{Pr cp1} Let $u_\e(x, t)$ be the unique global solution of
the CP for the regularized non-degenerate equation \eqref{cp4},
\ef{cp51}. Then, for
$t\in[0,T]$,  there exists a positive constant $K>0$,
independent of $\e$ and $T$, such that the following is satisfied:
\begin{enumerate}
\item[(i)] $\int\limits_\O |\nabla u_\e(x,t)|^2 \leq K$,\,\,
 $\int\limits_\O  |u_\e(x,t)|^2 \leq K$;
\item[(ii)] $\int\limits_\O  u_\e(x,t) \leq K$;
\item[(iii)] $\left\|h_\e \right\|_{L^2(\O \times [0,T])} \leq K$, with
$h_\e := \phi_\e(u_\e) \nabla \D u_\e$.
\end{enumerate}
\end{proposition}

  \noi{\em Proof.}
Firstly, multiplying \eqref{cp4} by $\D u_\e$, integrating in $\O
\times [0,t]$ for any $t\in[0,T]$, and applying the formula of
integration by parts yields
\begin{equation}
\label{cp6}
 \tex{
    \frac{1}{2}\, \int\limits_\O |\nabla u_\e(x,t)|^2 +
    \int\limits_0^t \int\limits_\O \phi_\e(u)|\nabla \D u_\e|^2 =
    \frac{1}{2}\, \int\limits_\O |\nabla u_\e(x,0)|^2,
    }
\end{equation}
thanks to the boundary conditions \eqref{i3}. Note that
\begin{equation*}
 \begin{matrix}
      \int\limits_\O [|\nabla u_\e(x,t+h)|^2 - |\nabla u_\e(x,t)|^2]
 \ssk\ssk\\
      =\,
      -\int\limits_\O [\D u_\e(x,t+h)+ \D u_\e(x,t)][u_\e(x,t+h)-u_\e(x,t)].
      \end{matrix}
\end{equation*}
Then, dividing that equality by $h$, passing to the limit as
$h\downarrow 0$, and integrating between 0 and any $t\in [0,T]$,
we find that\footnote{In fact, this is true from the beginning for
classic $C^\infty$-smooth solutions of (\ref{cp4}), but we will
need those manipulations in what follows.}
\begin{equation*}
 \tex{
    \int\limits_0^t \int\limits_\O \D u_\e u_{\e,t} =\frac{1}{2}\, \int\limits_\O |\nabla u_\e(x,t)|^2 -
    \frac{1}{2}\, \int\limits_\O |\nabla u_\e(x,0)|^2,
    }
\end{equation*}
which provides us with the necessary expression to obtain
\eqref{cp6}. Thus, from \eqref{cp6}, we have
\begin{equation}
\label{cp33}
 \tex{
    \int\limits_\O |\nabla u_\e(x,t)|^2 \leq K
    \quad \hbox{and} \quad
    \int\limits_0^t \int\limits_\O \phi_\e(u)|\nabla \D u_\e|^2 \leq
    K,
    }
\end{equation}
since both terms of the left-hand side in \eqref{cp6} are always
positive and the right-hand side is bounded by \eqref{i4}, for
a positive constant $K>0$ that is independent of $\e$. Then,
by Poincar\'e's inequality (when $\O$ is a bounded domain),
$$
 \tex{
   \int\limits_\O |u_\e(\cdot,t)|^2 \leq K \quad \mbox{for all} \quad t \in(0,T).
   }
   $$
   When $\O= \re^N$,
 for compactly supported data, all the solutions $u_\e(x,t)$ have
exponential decay at infinity, which makes the integrations
properly justified. Then the last inequality remains valid in
  certain $L^2_\rho(\re^N)$ and $H^1_\rho(\re^N)$ weighted spaces
    for an appropriately exponentially decaying weight $\rho$ (say $\rho(y)={\rm e}^{-a|y|^{4/3}}$).

On the other hand, in view of the finite propagation of
perturbations in the TFE--4, on a given interval $t \in [0,T]$, we
actually always deal with uniformly bounded supports of $u(x,t)$,
while $u_\e(x,t)$ exhibit just extra exponentially small tails at
infinity that do not affect the estimate.

Furthermore, from the mass conservation, we can also assure that
$$
 \tex{
\int\limits_\O u_\e(\cdot,t) \leq K \quad \mbox{for all} \quad t
\in(0,T).
 }
$$
Now, using \eqref{cp6} for $u_\e(\cdot,t) \in L^1(\O)$, we prove
the following estimate:
\begin{equation}
\label{cp12}
    \left\|h_\e \right\|_{L^2(\O \times [0,T])} \leq K,
    \quad \hbox{where} \quad h_\e = \phi_\e(u_\e) \nabla \D u_\e.
\end{equation}
To this end, from \eqref{cp33}, we find
\begin{equation}
\label{cp52}
 \begin{matrix}
 \tex{
    \int\limits_0^t \int\limits_\O
    [\e^n+(1-\e)(\e^2+u^2)^{\frac n2}]|\nabla \D u_\e|^2 \leq
    K, \quad \mbox{so that}
    }
    \ssk\ssk\\
 \tex{
    \e^n\int\limits_0^t\int\limits_\O |\nabla \D u_\e|^2 \leq K
    \quad \hbox{and} \quad
    \int\limits_0^t \int\limits_\O
    (\e^2+u^2)^{\frac n2}|\nabla \D u_\e|^2 \leq
    K,
    }
     \end{matrix}
\end{equation}
\noi since $\e \in (0,1)$, with a constant $K>0$ independent of
$\e$. Now, using H\"older's inequality,
\begin{equation*}
 \tex{
    \int\limits_0^t \int\limits_\O |h_\e|^2 \leq
    2\e^{2n} \int\limits_0^t \int\limits_\O
    |\nabla \D u_\e|^2 +2 \int\limits_0^t \int\limits_\O
    (\e^2+u^2)^{\frac n2}(\e^2+u^2)^{\frac n2}|\nabla \D u_\e|^2,
    }
\end{equation*}
and by \eqref{cp52}
(note also that  $u_\e(\cdot,t)
\in L^1(\O)$), we obtain \eqref{cp12}.
$\qed$

\ssk

Additionally, we get uniformly bounded estimates in $L^\iy$ for
the solutions of the TFE--4 equation \eqref{i1} as well as for the
solutions of the non-degenerate equation \eqref{cp4}.

Indeed, without loss of generality, we consider a solution
$u(x,t)$ of \ef{i1}, and, by scaling techniques (see
\cite{GMPKS,GW2}), prove its uniform $L^\iy$ {\em a priori} bound.

The function $u(x,t)$ is assumed to satisfy the equation in the
classic sense in the positivity and negativity domain, where it is
$C^\iy$ in $x$ (or may be even  analytic in $x$) by classic
parabolic theory \cite{Fr}, so the proof does not use particular
homogeneous smooth and sufficiently regular free boundary
conditions at the zero set. 

Implicitly, however, we have to assume
that such free boundary conditions cannot lead to some ``strong
singularities" on the free boundaries, which are not small in
$L^\iy$ (a kind of ``blow-up at the interfaces") and, hence, that
may affect standard interior parabolic regularity phenomena in the
uniform positivity/negativity domains. 

In other words, we assume
that, while solutions remain bounded and sufficiently smooth,
homogeneous free boundary conditions imposed (if any; e.g.,  no
conditions are assumed for the Cauchy problem in $\ren \times
\re_+$) move positive and negative humps of solutions sufficiently
smoothly and without strong ``collisions" and
``self-focusing-like" behaviour leading to $L^\iy$ blow-up at
interface points.

 The same analysis also directly applies
to the regularized uniformly parabolic problems like \ef{cp4},
where the classic solution is unique, so it can be used in
Theorem \ref{principal} and in other related results.

\begin{proposition}
\label{Pr.Glob} Under the above hypothesis,
 any solution of \ef{i1} is uniformly bounded.
    \end{proposition}

   \noi{\em Proof.} We argue by contradiction. Let there exist a
   monotone sequence $\{t_k\} \to +\iy$ (or $\{t_k\} \to T< \iy$,
   a finite blow-up time that makes no difference) and $\{x_k\} \subset \ren$ such that
    \be
    \label{gl.2}
     \tex{
    \sup_{(x,t) \in \ren\times (0,t_k)}|u(x,t)|= |u(x_k,t_k)|=C_k \to +\iy \quad
    \mbox{monotonically}.
    }
    \ee
    Then we rescale the TFE--4 and introduce a sequence
   of solutions $\{v_k(y,s)\}$ of the TFE--4:
    \be
    \label{gl.3}
    \begin{aligned}
     &
 \tex{
     u_k(x,t) \equiv u(x+x_k,t+t_k)=C_kv_k(y,s), \quad x=a_k y, \quad t=b_k
     s, \quad%%\hbox{with}\quad 
     \frac {b_k C_k^n}{a_k^4}=1}
     \\ \LongA
    &
 \tex{
     v_s=-\n \cdot (|v|^n \n \D v), \quad s>-\frac{t_k}{b_k} \to -\iy,
      \quad v_{k0}(y)= \frac 1{C_k}\,u_0(x_k+a_k y).
 }
    \end{aligned}
    \ee
    With this rescaling we are just performing a zoom around the point $(x_k,t_k)$, 
    in the region 
    $$B_\d(0)\times (-\t_k,0), \quad\hbox{with}\quad \d>0\quad \hbox{sufficiently small,}$$ 
    defining the problem in a ball 
    of radius $\frac{\d}{a_k}$ for the variable $y$ and in the interval 
    $$s\in \left(-\frac{t_k}{b_k},0\right).$$ 
    Note that we have a negative interval since 
    we have chosen $t+t_k$ and we approach the time from the left. Hence, the limiting problem will be defined in $\re^N\times  (-\infty,0)$. 
    The only assumption is that  the positive sequences satisfy 
    $$\{a_k\}  \to 0\quad \hbox{and} \quad \{b_k\} \to 0,$$ 
    for  scaling reasons.
The sequence of solutions $\{v_k(y,s)\}$ thanks to \eqref{gl.2} then satisfies that  \be
 \label{gl.4}
 |v_k(y,s)| \le 1 \andA |v_k(0,0)|=1, \quad\hbox{for all}\quad s \in \left[- \frac {t_k}{b_k},0\right).
 \ee
 Moreover, by the uniform estimate as in Proposition \ref{Pr
 cp1}(i), we have
  \be
  \label{gl.5}
   \tex{
 \int |\n_x u|^2 {\mathrm d}x \equiv \frac {C_k^2}{a_k^{N+2}}\,  \int |\n_y v_k|^2 {\mathrm d}y
 \le K \LongA \int |\n_y v_k|^2 {\mathrm d}y \le K
 \frac{a_k^{N+2}}{C_k^2}.
  }
  \ee
  Thus,  passing to limit, as $k \to +\iy$, along a subsequence, we are then supposed to have an
    {\em ancient} solution $v_k \to v(y,s)$ defined for all $s<0$.
According to the last estimate in \ef{gl.5}, we will get a zero
solution in the limit provided that
 \be
 \label{gl.6}
  \tex{
 \kappa_k \equiv \frac {a_k^{N+2}}{C_k^2} \to 0.
}
 \ee
 Thus, we have to satisfy the following two properties only:
  \be
  \label{gl.7}
   \left\{
   \begin{aligned}
 & b_k=\kappa_k^4 C_k^{\frac 8{N+2}-n} \to 0, \\
 & a_k=\kappa_k^{\frac 1{N+2}} C_k^{\frac 2{N+2}} \to 0.
 \end{aligned}
 \right.
  \ee
Evidently,
 both can be done easily provided that $\kappa_k$ decays
 sufficiently fast.
Hence, since passing to the limit, as $k \to +\iy$, in \eqref{gl.3} for bounded solutions in the
uniform positivity (negativity) domains of a parabolic equation is
straightforward, the ancient solution is $v \equiv 0$. This
contradicts the last assumption in \ef{gl.4}, implying, by the
interior parabolic regularity, that $v(y,0)$ must be non-trivial
in a neighbourhood of $y=0$. This simply means that the TFE--4
does not have an internal mechanism to support indefinite growth
(or blow-up) of solutions.
 $\qed$
 
 \vspace{0.2cm}
 
 {\bf Remark.} \rm{This type of scaling proof is rather standard, based on the famous Gidas--Spruck blow-up Method \cite{GS81} for elliptic problems. 
 For parabolic problems see \cite{GMPKS} where a Kuramoto--Sivashinsky equation is analysed and \cite{GW2} for higher-order equations with absorption terms.}

%%%%%%%%%%%%%%%%%%%%%%%%%%%%%%%%%%%%%%%%%%
  \subsection{Passing to
the limit}

To conclude this section, we show the existence of a weak solution of the CP
for the degenerate parabolic TFE--4 \eqref{i1} by passing to the
limit: 
\begin{enumerate}
\item[(i)] as $\e$ goes to zero, 
\end{enumerate}
and also, 
\begin{enumerate}
\item[(ii)] as  $n \to 0^+$.
\end{enumerate}

\ssk

  %%%%%%%%%%%%%%%%%%%%%%%%%%%%%%%%%%%%%%%%%%%%
\noi\underline{\sc Passing to the limit $\e \to 0^+$}. By
Proposition\;\ref{Pr cp1}, for bounded supports $\O$, and the
Aubin--Lions Lemma, the embedding 
$$H_0^1(\O\times (0,T)) \hookrightarrow L^2(\O\times (0,T)),$$ 
is compact. Then, we can extract a convergent
subsequence in $L^2(\O\times (0,T))$ as $\e\downarrow 0^+$ of
solutions of \eqref{cp4}, \ef{cp51}, labelled again by
$u_\e(x,t)$, so
\begin{equation}
\label{cp11}
    \lim_{\e \rightarrow 0^+} \big\| u_\e(\cdot,t) - U(\cdot,t)\big\|_{L^2(\O\times (0,T))}
    =0,
\end{equation}
where, for convenience, unlike in \ef{Fem1}, this CP-solution of
the TFE--4 \ef{i1} is denoted by the capital $U(\cdot,t)$.

Thus,
the convergence of the regularized solutions of the CP
\eqref{cp4}, \ef{cp51} is strong in $L^2(\O\times (0,T))$. In
$\re^N$, we will use the appropriate $L^2_\rho$ and Sobolev spaces
$H^1_\rho$ when necessary.

Note again that one difficulty
we face is whether this limit depends on the taken subsequence or
not. 

In other words, this analysis, just with the limit $\e \to 0^+$, does not include
 any uniqueness result, which is expected to be a more difficult
 open problem for such nonlinear degenerate parabolic TFEs in
 non-fully divergence form and with non-monotone operators.

 However, the principal issue of the analytic regularization via
 \eqref{cp4}, \ef{cp51} is that it is expected to lead to smoother CP-solutions at the
 interface than those for the standard FBP, at least for $n < \frac 32$. 
 
 The difference is that
using the analytic regularized family $\{u_\e\}$ and imposing
 \eqref{i3} it is assumed to guarantee that on the interface
 (being sufficiently smooth for simplicity),
  \be
  \label{int4}
  \tex{
   \frac{\p^2 u}{\p {\bf n}^2} =0 \quad \hbox{a.e., at least, for $n \in (0,1]$.}
   }
    \ee
In fact, proper oscillatory solutions of the CP are assumed to
exhibit even more regularity at smooth interfaces \cite{EGK2} (cf.
\ef{Conv2}): again in the sense of non-oscillatory
envelopes\footnote{Or in the sense of partial limits along
subsequences staying ``sufficiently away" from less smooth
transversal zeros concentrating at the interfaces.},
 \be
 \label{int5}
\tex{
   \frac{\p^l u}{\p {\bf n}^l} =0, \quad \mbox{where} \quad
   l=\big[\frac 3n \big]-1.
   }
    \ee
Therefore, as $n \to 0^+$, the smoothness of such solutions at the
interfaces increases without bounds. Obviously, this is not the
case for the FBP (a ``positive obstacle" one) with standard
conditions \eqref{i3} and a usual quadratic (``parabolic") decay
at the interfaces for any $n < \frac 32$, \cite{Gia08}.

\begin{lemma}
Let $u_\e(x, t)$ be the unique global solution of
the CP for the regularized non-degenerate equation \eqref{cp4},
\ef{cp51}. Then,
$$ \lim_{\e \rightarrow 0^+} \big\| u_\e(\cdot,t) - U(\cdot,t)\big\|_{L^2(\O\times (0,T))}
    =0,
    $$
 with $U(\cdot,t)$ a weak solution of the TFE--4 \eqref{i1} satisfying
$$
 \tex{    \int\limits_0^T \int\limits_\O \var_t U
    +\int\limits_0^T \int\limits_\O \nabla\var \cdot (|U|^n \nabla \D U)
    =0 \quad \mbox{for any} \quad \var \in C_0^\infty(\bar\O \times [0,T]).
 }
$$
\end{lemma}

 \noi{\em Proof.}
Thus, as above and customary, multiplying \eqref{cp4} by a test
function $\var \in C_0^\infty(\bar\O \times [0,T])$ and
integrating by parts in $\O \times [0,T]$
   gives
\begin{equation*}
 \tex{
    -\int\limits_0^T \int\limits_\O \var_t u_\e -
    \int\limits_0^T \int\limits_\O \nabla\var \cdot (\phi_\e(u)\nabla \D u_\e)
    =0.
    }
\end{equation*}
 Next, dealing with this equality and  assuming $\phi_\e(u)$ of the form \eqref{cp51}, we find
\begin{equation}
\label{cp13}
 \tex{
    \int\limits_0^T \int\limits_\O \var_t u_\e + \e^n \int\limits_0^T \int\limits_\O
    \nabla\var \cdot \nabla \D u_\e
    +(1-\e)\int\limits_0^T \int\limits_\O \nabla\var \cdot (
    (\e^2+u^2)^{\frac n2}\nabla \D u_\e) =0.
    }
\end{equation}
Applying H\"{o}lder's inequality, it is clear from \eqref{cp6}
that there exists a subsequence labelled by $\{\e_k\}$, so that
the second term of \eqref{cp13} approximates zero as $\e_k
\downarrow 0^+$ for small $n>0$,
 \be
 \label{cp13N}
 \tex{
    \Big|\e_k^n \int\limits_0^T \int\limits_\O \nabla\var \cdot \nabla \D u_{\e_k}\Big|
    }
 \tex{
    \leq \e_k\Big(\e_k^{2(n-1)} \int\limits_0^T \int\limits_\O |\nabla \D
     u_{\e_k}|^2\Big)^{\frac 12}
    \Big( \int\limits_0^T \int\limits_\O |\nabla\var|^2\Big)^{\frac 12}
 }
    \\
     \leq K \e_k \downarrow 0,
 \ee
as $\e_k \downarrow 0^+$, for a positive constant $K>0$.

Moreover, on the ``good subsets" of uniform non-degeneracy
  \be
  \label{G11}
  \mathcal{G}_{\e,\d}:=\{(x,t)\in \O\times
[0,T]\,;\, |u_\e(x,t)|> \d>0\}
  \ee
   for any fixed arbitrarily small $\d>0$,
it is clear that the limiting solution as $\e \to 0^+$ is a weak
solution of the TFE--4 \eqref{i1}.

Moreover, this represents a classical interior regularity result for
parabolic equations on the sets of their uniform parabolicity.
Indeed, by the regularity for the uniformly parabolic equation
\eqref{cp4},
we obtain that $u_{\e,t}$, $\nabla u_{\e}$,
and $\phi_\e(u) \nabla \D u_{\e,x}$
converge on
compact subsets of $\mc{G}= \mc{G}_{0,0}$ (see below), at least in a weak sense. 

In general, it is not that difficult
to see that, as $\e=\e_k \to 0^+$ (along the lines of classic
results in \cite{BF1} and a number of later related others), we
obtain a weak solution of the TFE--4, i.e.,
\begin{equation}
\label{cp14}
 \tex{
    \int\limits_0^T \int\limits_{\mathcal{G}} \var_t U
    +\int\limits_0^T \int\limits_{\mathcal{G}} \nabla\var \cdot (|U|^n \nabla \D U) =0,
 }
\end{equation}
where $U(x,t)$ is the limit obtained through \eqref{cp11},
${\mathcal G}={\mathcal G}_{0,0}$ is associated with $U(x,t)$
(recall that it is  smooth away from the nodal set by the
interior parabolic regularity), and naturally assuming that $\var
\in C^\iy_0({\mathcal G})$.

\ssk

However, in the ``bad subsets"
  \be
  \label{B11}
   \mathcal {B}_{\e,\d}:=\{(x,t): \; |u_\e(x,t)|\leq \d\} \quad \hbox{for
 sufficiently small} \quad  \d \geq 0,
   \ee
   we must take $\e>0$ also
 small enough  and depending on $\d$.
Indeed, let us fix an $\e$ such that $0<\e\leq \d$. Thus, applying
H\"{o}lder's inequality to the third term in \eqref{cp13} over the
set \ef{B11}, where $|u_\e|\leq \d$, we have that
\begin{align*}
 \tex{
    \Big|  \int\limits_0^T \int\limits_{\{|u_\e|\leq \d\}}
 }
     &
      \tex{
       \nabla\var \cdot
    ( (1-\e) (\e^2+u_\e^2)^{\frac n2} \nabla \D u_\e)
    \Big|
    }
    \\ &
 \tex{
    \leq \Big(\int\limits_0^T \int\limits_{\{|u_\e|\leq \d\}} |\nabla\var|^2 \Big)^{\frac 12}
     \Big(\int\limits_0^T \int\limits_{\{|u_\e|\leq \d\}} (1-\e)^2 (\e^2+u_\e^2)^n
     |\nabla \D u_\e|^2\Big)^{\frac 12}.
     }
\end{align*}
Then, since $\var \in C_0^\infty(\bar\O \times (0,\infty))$ and
$\e\in (0,1)$, we get
\begin{equation*}
 \tex{
    \Big|  \int\limits_0^T \int\limits_{\{|u_\e|\leq \d\}}
 }
      \tex{
      \nabla\var \cdot
    ( (1-\e) (\e^2+u_\e^2)^{\frac n2} \nabla \D u_\e) \Big|
 }
 \tex{
  \leq   C \Big(\int\limits_0^T \int\limits_{\{|u_\e|\leq \d\}} (1-\e) (\e^2+u_\e^2)^n
     |\nabla \D u_\e|^2\Big)^{\frac 12}
     }
\end{equation*}
for a positive constant $C>0$. Using that $|u_\e|\leq \d$, by
\eqref{cp33},  we find that
\begin{align*}
 \tex{
    \Big|  \int\limits_0^T \int\limits_{\{|u_\e|\leq \d\}}
 }
    &
 \tex{
    \nabla\var \cdot
    ( (1-\e) (\e^2+u_\e^2)^{\frac n2} \nabla \D u_\e) \Big|
 }
    \\ &
     \tex{
   \leq  C \Big(\int\limits_0^T \int\limits_{\{|u_\e|\leq \d\}} (1-\e) (\e^2+\d^2)^{\frac n2}
     (\e^2+u_\e^2)^{\frac n2}     |\nabla \D u_\e|^2\Big)^{\frac 12},
      }
\end{align*}
and, hence, using \eqref{cp52},
\begin{equation}
 \label{nn11}
 \tex{
    \Big|  \int\limits_0^T \int\limits_{\{|u_\e|\leq \d\}} \nabla\var \cdot
    ( (1-\e) (\e^2+u_\e^2)^{\frac n2} \nabla \D u_\e) \Big|\leq
     C_1 \d^{\frac n2}
       \to 0,
     }
\end{equation}
for a constant $C_1>0$, provided that $\d \to 0$ as $\e \to 0^+$.
This means that integration over ``bad subsets" \ef{B11} leaves no trace in the
final limit. So, the integral identity \ef{cp14} for such a
CP-solution $U$  holds for arbitrary domains ${\mathcal G}$, so
that $U$ is a true weak solution of the CP \ef{i1}.
$\qed$

\ssk

 {\bf Remark.} \rm{Although, we are able to prove the existence of a limit in a weak sense and 
 the continuous deformation from the parabolic bi-harmonic equation \eqref{s1} to the TFE--4 \eqref{i1} through 
 the regularized non-degenerate equation \eqref{cp4}, \ef{cp51}, this also proves that a  standard integral identity is not enough to specify directly the solutions of the CP.
  Indeed, due to free boundary conditions \ef{i3}, the weak formulation \ef{cp14} holds for any FBP-solution.}

\vspace{0.4cm}

\ssk

%%%%%%%%%%%%%%%%%%%%%%%%%%%%%%%%%%%%%%%%%%%%%%%%%%%%%%%%%%%%%%%%
 \noi\underline{\sc Passage to the limit $n \to 0^+$}. In addition, the
estimate \eqref{nn11} shows the actual rate of convergence if we perform a
second ``homotopy"
 limit as $n \to 0$,  together with $\e \to 0$, in the analytic
 approximating flow \eqref{cp4}, \ef{cp51} in order to obtain, in this limit over \eqref{cp13}, weak
 solutions of the bi-harmonic equation \eqref{s1} and,
 hence, classical solutions by standard parabolic theory. For instance, to get such a convergence, one needs
  \be
  \label{nn12}
   \tex{
    \mbox{for} \quad \d \sim \e, \quad
  n=n(\e) \to 0^+ \,\,\, \mbox{such that} \quad \e^{\frac{n(\e)}2}
  \to 0^+.
   }
   \ee
Indeed, setting
\begin{equation*}
    \lim_{\e \rightarrow 0^+} \e^{\frac{n(\e)}2} = 0, \quad \hbox{then} \quad
     \lim_{\e \rightarrow 0^+} n(\e) \ln \e= -\iy.
\end{equation*}
 Then, we need the following:
\begin{equation*}
    \tex{ n(\e) \ln \e \ll -1 \LongA n(\e) \gg \frac{1}{|\ln \e|},}
\end{equation*}
which will provide us with the limit $n(\e) |\ln \e| \to + \iy$,
and convergence of solutions, at least, in a weak sense.

However, this is not the end of the problem: indeed, under the
condition \eqref{nn12} on the parameters, we definitely arrive, in
the limit as $\e,\, n(\e) \to 0^+$, at a weak solution of the
bi-harmonic equation \ef{s1} written in the following ``mild"
form:
\begin{equation}
\label{cp14N}
 \tex{
    \int
    \limits
    _0^T \int
    \limits
    _\O \var_t U
    +\int
    \limits
    _0^T \int
    \limits
    _\O \nabla\var \cdot \nabla \D U
    =0.
 }
\end{equation}
This is not a full (strongest) definition of weak solutions, since
it assumes just a single integration by parts, so formally allows
us to obtain, in the limit, a positive solution of the ``obstacle" FBP
for \eqref{s1} with the corresponding conditions \eqref{i3} (with
$n=0$), which can be constructed by a ``singular" regularization as
in \cite{BF1}. Other solutions of different FBPs cannot be ruled
out either.

Thus,  it turns out that such a more or less standard passage to
the limits in the integral identities {\em in principal cannot}
recognize the desired difference between oscillatory solutions of
the CP and others (possibly, non-negative or even positive ones)
of the FBP posed for the TFE--4  \eqref{i1}. This emphasizes in
general how hard a proper definition of the CP-setting can be.
 Therefore, a
stronger version
 of homotopic limits is crucially necessary towards a proper
 identification of the CP-solutions.

Nevertheless, this first preliminary step in homotopy analysis
declares some useful estimates and bounds on the parameters of
regularization such as \eqref{nn12}, which are  necessary for a
correct passing to the limit to get sign changing solutions of the
linear bi-harmonic flow. 

Moreover, even now, the homotopy concept
as a connection to the linear PDE \eqref{cp4} can  describe the
origin (at $n=0$) of the oscillatory solutions of TFEs and, hence,
establish a transition to the {\em maximal regularity
solutions}  of the CP \eqref{i1}. Inevitably, bearing in mind the
oscillatory character of the kernel $F(|y|)$ of the fundamental
solution \ef{s3}, the proper solutions of the CP are going to be
oscillatory near finite interfaces, at least, for all small $n >0$.

%%%%%%%%%%%%%%%%%%%%%%%%%%%%%
\section{Branching through  homotopy \ef{phi2} as $n \to 0^+$ and $\e \to 0^+$}
 \label{S3}

 %%%%%%%%%%%%%%%%%%%%%%%%%%%%%%%%%%%%%%%%%

Now, performing a double limit as $n \to 0^+$ and $\e \to 0^+$, we
ascertain the well-posedness of the TFE--4
\eqref{i1} via a homotopy deformation from the solutions of the {\em bi-harmonic
equation} \eqref{s1}, which are known to be  oscillatory, to the solutions of the TFE--4 \eqref{i1}. This will
provide us with one of the main results of this paper.

 As a key example, consider the regularized CP \ef{cp4}, \ef{phi2}, with the same
 data:
  \be
  \label{3.1}
u_\e: \quad  u_t=- \n \cdot ((\e^2+u^2)^{\frac n2} \n \D u), \quad
u(x,0)=u_0(x),
 \ee
 where smooth bounded data $u_0(x)$ are compactly supported, or, at
 least, has exponential decay at infinity.
 Then, quite naturally, due to local parabolic properties, we may assume that, on any bounded fixed
 interval $t \in [0,T]$ and for any small $\e>0$,
\be
 \label{3.11}
 u_\e(x,t) \quad \mbox{has exponential decay as $x \to \iy$}.
 \ee
 Moreover, we may also suppose that $u_\e(x,t)$ has a natural spatial form
 of a ``hump" concentrated near the origin on the given fixed time-interval $t \in [0,T]$ (in other words, for
 such a uniformly parabolic problem, there is no fast propagation
 at $x=\iy$, since parabolic flows do not allow this by
 interior regularity theory).

 Contrary to what we did above we now choose 
 $$n \to 0^+,$$ 
 as the main deformation parameter, and,
 for simplicity, will properly choose
 $$\e=\e(n) \to 0.$$
  Thus, we perform a passing to a single limit, though similar
 estimates can be obtained within a double (or even triple if necessary) limit strategy.

To this end, we  rewrite this PDE in a perturbed version,
  \be
  \label{3.1.1}
  u_t= - \D^2 u + \n \cdot \bigl([1-(\e^2+u^2)^{\frac n2}] \n \D u \bigr),
 \ee
 and next, we use its equivalent integral form:
  \be
  \label{3.3}
   \begin{matrix}
  u_\e: \quad u(t)= b(t) * u_0 + \int\limits_0^t \n b(t-s)* F_{n,\e}(u(s)) \n \D
  u(s)\, {\mathrm d} s,
   \ssk\ssk\\
   \mbox{where} \quad F_{n,\e}(u):= 1-(\e^2+u^2)^{\frac
  n2}.
  \end{matrix}
   \ee
 The convergence to the well-posed bi-harmonic solutions
 will crucially depend on the weak limit in the second perturbation term
 in \ef{3.3}:
  \be
  \label{3.4}
 F_{n,\e}(u) \equiv 1-(\e^2+u^2)^{\frac
  n2}\rightharpoonup 0 \asA n, \, \e(n) \to 0^+.
   \ee
 In other words, we need to clarify conditions, under which
 there occurs branching of proper CP-solutions of the TFE--4 at
 $n=0^+$ from regular bounded classical solutions of the linear
 problem \ef{s1}. Recall that this requires branching theory for a
 quasilinear degenerate {\em partial} differential equation, where
 a good functional setting will be key.

 \vspace{0.4cm}

  \noi\underline{\sc Key assumption}. To perform our branching analysis, one needs to verify
 the following expansion:
  \be
  \label{3.5}
   \tex{
   F_{n,\e}(u) \equiv 1-(\e^2+u^2)^{\frac
  n2}= - \frac n2 \, \ln (\e^2+u^2)(1+o(1)) \asA n \to 0^+
  }
  \ee
 on a fixed family $\{u_\e(x,t)\}$ of uniformly bounded and smooth
 solutions. Checking \ef{3.5} on bad sets \ef{B11}, i.e., for $u
 \approx 0$, yields the demand (cf. \ef{nn12} if $\e \ll {\mathrm e}^{- \frac 1n}$)
  \be
  \label{3.6}
   \fbox{$
  n \, |\ln \e(n)| \to 0 \asA n \to 0^+.
   $}
  \ee
  This is the key assumption on the regularization parameter $\e=\e(n)>0$ and its relation with the main parameter $n$.

  \ssk

  %%%%%%%%%%%%%%%%%%%%%%%%%%%%%%%%%%%%%%

  \subsection{Proof of Theorem\;\ref{principal}}

  Under condition \ef{3.6}, we perform a rather standard branching analysis which mimics that in the elliptic case shown in
   \cite{TFE4PabloVictor}. Actually, the integral equation \ef{3.3}
   provides a nice opportunity to get a direct $n$-expansion of the
   solution to guarantee branching at $n=0^+$.

Hence, substituting \eqref{3.5} into \eqref{3.3}, we have that
\be
\label{lya}
    \tex{ u(t)= b(t) * u_0 -\frac{n}{2} \int\limits_0^t \n b(t-s)*
    \ln(\e^2+u^2) \n \D u(s)\, {\mathrm d} s + o(n) \asA n \to 0^+.}
\ee
Next, in order to derive the actual rate of convergence to
ascertain such a branching behaviour for small $n>0$, we assume the
following splitting expression:
\be
\label{spl}
    \tex{u = \tilde{u} + V,}
\ee
where $\tilde{u}$ is the unique solution \eqref{s2} of the
parabolic bi-harmonic equation \eqref{s1} with the same initial
data and with $V$ being a small perturbation. Thus, setting
 \be
 \label{spl11}
\tex{ V:= n \varphi_1 + o(n),}
 \ee
  and substituting
\eqref{spl} into \eqref{lya}, omitting the terms $o(n)$ when
necessary, yields
\begin{align*}
    \tex{ \tilde{u}(t) + n \varphi_1(t)} &  = b(t) * u_0
     \\
     & \tex{-\frac{n}{2} \int\limits_0^t \n b(t-s)*
    \ln\big(\e^2+\tilde{u}^2(s) + 2 n \tilde{u}(s) \varphi_1(s) + n^2 \varphi_1^2(s)\big) \n \D \tilde{u}(s)\, {\mathrm d} s }
 \\ &
    \tex{-\frac{n^2}{2} \int\limits_0^t \n b(t-s)*
    \ln\big(\e^2+\tilde{u}^2(s) + 2 n \tilde{u}(s) \varphi_1(s) + n^2 \varphi_1^2(s)\big) \n \D \varphi_1(s) \, {\mathrm d} s.}
\end{align*}
Furthermore, by the convolution operation properties, we find that
 \[ 
\begin{split}
    \tex{ \tilde{u}(t) + n \varphi_1(t)} & \tex{  = b(t) * u_0 } \\ & \tex{ -\int\limits_0^t \n b(t-s)
    *\frac{n}{2} \ln(\e^2+\tilde{u}^2(s) + 2 n \tilde{u}(s) \varphi_1(s) + n^2 \varphi_1^2(s)) \n \D \tilde{u}(s)\, {\mathrm d} s. }
\end{split}
 \] 
 Overall, this gives the following corrections in the expansion \ef{spl},
 \ef{spl11}:
 \be
 \label{phi12}
  \tex{
   \varphi_1(t)= -  \int_0^t  \n b(t-s) *  \ln | \tilde u(s)|\n \D \tilde u(s) \, {\mathrm d} s.
   }
   \ee
Thus, the above asymptotic expansion technique must
   assume that
 the convolution in \ef{phi12} is always finite, e.g.,  that
 $$\ln | \tilde u(y,s)|
\in L^1_{\rm loc}(\ren),$$
for any $s>0$, i.e., $\tilde u(y,s)$ does
not have zeros with an exponential decay in a neighbourhood. In
particular, this is true if the solutions \ef{s2} have {\em
transversal} zeros a.e. (this always happens in the radial
geometry, for instance). The latter is rather plausible, but
difficult to prove for general solutions of the bi-harmonic
equation \ef{s1}, so we need to include such an assumption.

 Subsequently:
 \begin{enumerate}
 \item[(i)] We need to prove that the perturbation
in \ef{3.3} is asymptotically small. Obviously, this is guaranteed
by the uniform estimate in Proposition \ref{Pr cp1} (iii) in the
domain $\{|u| \ge z_*\}$, where $z_*$ is such that
 \be
 \label{ln1}
  |\ln z| \le z^{\frac n2} \quad \mbox{for} \quad z \ge z_*.
  \ee

Consider now the double integral in \ef{3.3} over the domain,
where $$\{\e^2 \le \e^2+u^2 \le z_*\}.$$ Then, the maximal
singularity of the term $\ln(\e^2+u^2)$  therein is achieved at
$u=0$ and is $\sim \ln \e$. Concerning the third derivative
therein, quite similarly, as in the above case, we obtain from the
same estimate of the Proposition \ref{Pr cp1}(iii) that
 \be
 \label{est77}
  \tex{
  \iint|\n \D u_\e|^2 \le K\e^{-n},
   }
   \ee
   i.e., overall, the second perturbation term in \ef{3.3} has the
   order, at most,
    \be
    \label{est78}
    O(n \e^{-n}) \to 0 \asA n \to 0^+.
    \ee
Indeed, by \ef{3.6}, $\e^{-n}(n) \to 1$, so that \ef{est78} holds.

\ssk

\item[(ii)] As we have shown, the improved expansion \ef{spl11} requires
convergence of \ef{phi12}, as $n\to 0^+$, which is difficult to verify for
sufficiently arbitrary solutions of the bi-harmonic equation. 

It
is worth mentioning that, thanks to the $x$- and $t$-analyticity of these
solutions for $t \in (0,T]$, divergence of such integrals (if
any), i.e., violation of such a rate of divergence due to formation
of ``flat" multiple zeros of $\tilde u(x,t)$, can occur at a
finite number of points $(x,t) \in \ren \times (0,T]$, so that
this expansion is expected to hold, at least, in the a.e. sense
for sufficiently arbitrary initial data. $\qed$
\end{enumerate}

 \ssk

 \noi\underline{\sc Conclusion}.
Simultaneously, the expansion \ef{spl}, \ef{spl11} presents a
proper definition of the CP-solutions of the TFE--4 \ef{i1}:
\begin{align*}
& \hbox{ {\em solutions of the Cauchy problem for \ef{i1} are
those, which can be deformed} } \\ &  \hbox{\em{ as $n \to 0^+$
via the analytic path \ef{3.1}, \ef{3.6} to the unique solutions
\ef{s2}}} \\ &  \hbox{\em{ of the bi-harmonic equation \ef{s1}
with the same data.}}
\end{align*}
Inevitably, for small $n>0$,
such correct and well defined CP-solutions must inherit oscillatory and other
sign-changing properties of this linear flow.

In other words, according to a homotopy approach, a proper setting
of the CP for the TFE--4 \ef{i1} requires the whole set of
solutions 
$$\{u(x,t), \, n>0\},$$ 
or, hence, a full 2D set
$$\{u_\e(x,t), \, \e>0,\, n>0\},$$ 
of regularized solutions.
Therefore, a successful definition of the CP for an {\em
individual} TFE--4 for a {\em fixed} value $n>0$ may not be
efficient or available.

%%%%%%%%%%%%%%%%%%%%%%%%%%%%%%%%%%%%
 \section{Towards uniqueness: boundary layer and Riemann's problems as $\e \to 0^+$}
 \label{S5}

 %%%%%%%%%%%%%%%%%%%%%%%%%%%%%%%%%%%%%%

We continue to study the regularized problem \ef{3.1}, and now we
concentrate on the crucial limit \ef{Fem1}. As we have mentioned,
a {\em conventional} (i.e., associated with the regularization
\ef{phi2} currently applied) uniqueness naturally requires that
 \be
 \label{5.1}
 \limsup u_\e(x,t)= \liminf u_\e(x,t) = u(x,t) \asA \e \to 0^+.
  \ee
Actually, this demands knowing that no essential $\e$-oscillations
of $u_\e(x,t)$ occur in the ``bad sets" \ef{B11}, which are close
to the nodal set of the limit solution $u(x,t)$. Indeed, if the
behaviour of $u_\e$ as $\e \to 0^+$ is non-oscillatory enough in
these singular sets, then such oscillations cannot occur in
principle in any good set of uniform positivity of solutions,
where a strong interior parabolic regularity is in charge.

Therefore, proving the uniqueness \ef{5.1} assumes a deeper
understanding, than before, of {\em boundary layers} that occur
close to nodal sets of the solutions, and mainly close to
interfaces, where the strongest {\em singularity via the equation
degeneracy} occurs. This leads to very difficult asymptotic
questions, which are naturally referred to as to {\em Riemann's
problems} for the TFE--4, each of which is associated with a
prescribed type of singularity at the interface. We now present a
short discussion of one such a Riemann problem.

%%%%%%%%%%%%%%%%%%%%%%%%%%%%%%
 \subsection{First (artificial) Riemann's problem}

 %%%%%%%%%%%%%%%%%%%%%%%%%%%%%%%%%%%%

 The regularized equation \ef{3.1}
has the clear advantage that $\e$ can be scaled out by a number of
scaling transformations:
 \be
 \label{5.2}
  \tex{
 u_\e(x,t)= \e v_\e(y,\t) \whereA y= \frac x{\e^\a}, \quad \t= \frac
 t{\e^{\b}},\,\,\, \b=4 \a -n,
  }
  \ee
  and $\a \in \re$ is an arbitrary exponent. Then $v=v_\e$ solves
  an $\e$-independent equation with data $v_{0\e}$ now depending on
  $\e$:
   \be
   \label{5.3}
    \tex{
   v_\e: \quad v_\t=- \n_y \cdot ((1+v^2)^{\frac n2} \n_y \D_y v),
   \quad v_{0\e}(y)= \frac 1 \e \, u_0(\e^\a y).
   }
   \ee
Consider a particular simpler case $\a=\frac n4$, for which $\b=0$
in \ef{5.2}, i.e., the time variable is not under scaling. Then we
have the following problem:
 \be
 \label{5.4}
  \tex{
  u_\e(x,t)= \e v_\e\bigl(\frac x{\e^{n/4}},t \bigr), \,\,\,
   v_t=- \n \cdot ((1+v^2)^{\frac n2} \n \D v),
   \quad v_{0\e}(y)= \frac 1 \e \, u_0(\e^{\frac n4} y).
    }
 \ee
 To formulate a particular Riemann problem, let us assume that
 $0 \in \partial \,{\rm supp}\, u_0$ and
  \be
  \label{5.5}
  \tex{
  u_0(x) = \chi\bigl(\frac x{|x|}\bigr) |x|^{\frac 4n}(1+o(1))
   \asA x \to 0,
   }
   \ee
   where $\chi(s) \not \equiv 0$ is a continuous function on the unit sphere. In other words, \ef{5.5}
   fixes a necessary type of degeneracy (vanishing) of initial data close to
   the origin being an interface point at the initial moment $t=0$.

   Later on, such a Riemann problem must be attached
to {\em every} point $x_0 \in \partial {\rm supp}\, u_0$ of the
initial interface, by changing $x \mapsto x-x_0$ in \ef{5.5}, or
with other data asymptotics, where other scaling transformation
from \ef{5.2} (with a different $\a$) can be applied.

 Thus, for data as in \ef{5.5}, after scaling as in \ef{5.4}, we obtain the following
 rescaled initial data: uniformly on compact subsets in $y$,
  \be
  \label{5.6}
   \tex{
  v_{0\e}(y)= \chi\bigl(\frac y{|y|}\bigr) |y|^{\frac 4n}(1+o(1)) \to
\chi\bigl(\frac y{|y|}\bigr) |y|^{\frac 4n} \equiv v_0(y)   \asA
\e \to 0,
 }
 \ee
  i.e., there exists a {locally} finite  limit  initial data. This gives {\em Riemann's problem}
 as a first approximation of a proper solution:
  \be
  \label{5.7}
   \tex{
   v_t=- \n \cdot ((1+v^2)^{\frac n2} \n \D v), \quad
    v_0(y)= \chi\bigl(\frac y{|y|}\bigr) |y|^{\frac 4n}.
    }
     \ee
 As usual, the first question is to check whether Riemann's problem \ef{5.7} is globally well-posed or 
 admits a singularity. It turns out that the latter holds:

 \begin{proposition}
 \label{P.R1}
 In general, Riemann's problem \ef{5.7} admits blow-up in finite time.
  \end{proposition}
 \noi{\em Proof.}
To see this, we perform some easy estimates. Namely, as is well
known, such a kind of blow-up is due to a fast growth of data in
\ef{5.7} as $y \to \iy$, where $v$ also gets arbitrarily large and
hence the PDE can be approximated by the initial TFE--4,
 \be
 \label{5.8}
 v_t=- \n \cdot(|v|^n \n \D v).
  \ee
Looking for separate variable solutions of \ef{5.8} yields
 \be
 \label{5.9}
 v(y,t)= \psi(t) \rho(y) \LongA \psi'=\psi^{n+1}, \quad \rho = - \n \cdot(|\rho|^n \n \D \rho).
  \ee
 The first ODE for $\psi>0$ yields a typical blow-up behaviour,
  \be
  \label{5.10}
   \psi(t)=[n(T-t)]^{-\frac 1n} \to \iy \asA t \to T^-,
    \ee
    where $T>0$ is a finite blow-up time. Dealing with the second nonlinear
    elliptic equation in \ef{5.9}, it is easy to see that such an
    equality assumes the following algebraic balance:
     \be
     \label{5.11}
     \rho(y) \sim |y|^{\frac 4n} \asA y \to \infty,
      \ee
      i.e., the same growth as data $v_0(y)$ in \ef{5.7}.
Moreover, the equation for $\rho$ admits an explicit radial
solution of this form,
 \be
 \label{5.12}
 \rho(y)=C_* |y|^{\frac 4n} \inB \ren \whereA C_*=C_*(n,N)={\rm
 const.}>0.
  \ee
Of course, some special changing sign behaviour of the angular
function $\chi(s)$ in data $v_0$ can result in global solvability
of Riemann's problem \ef{5.7}. However, since the data is assumed to
be arbitrary, we cannot rely on this. $\qed$

  \ssk

The blow-up result in Proposition \ref{P.R1} is precisely the sign
of possible {\em strong oscillations} that may occur as $\e \to 0$
close to such a degeneracy point. Note that, unlike second-order
heat equations, blow-up in higher-order parabolic equations
naturally means occurring infinitely many sign changes as $t \to
T^-$, i.e., solutions get infinitely oscillatory close to blow-up
time; see \cite{GS1S-V} and \cite{CG2m} as two examples.  In this
case, \ef{5.1} may become non-achievable.

 \ssk

On one hand, this discrepancy does not mean the actual blow-up of
the corresponding regularizing  family $\{u_\e\}$ as $\e \to 0$.
This example for special (and somehow artificial) data \ef{5.5}
just shows that such an independent asymptotic analysis, such as via
\ef{5.7}, of the boundary layer is not allowed, as well as a full study
of the $\e$-dependent problem \ef{5.4} (including a  careful
control of the behaviour at $y \to \iy$ as $\e \to 0$, to avoid
blow-up and essential oscillations).

On the other hand, we cannot exclude a possibility of the actual
existence of non-small $\e$-oscillations of $u_\e(x,t)$ in a
neighbourhood of the interface point for data \ef{5.5} {\em
provided that such a behaviour remains for some, possibly,
extremely small time intervals}. However, and fortunately, this is
not the case, and the generic behaviour near interfaces has
different scaling laws. Then, the above blow-up solutions just
describe a fast transition from data \ef{5.5} to a solution
$u(x,t)$ for arbitrarily small $t>0$, which already gains the
generic oscillatory interface behaviour to be described next.

%%%%%%%%%%%%%%%%%%%%%%%%%%%%%%%%%%%%%%%%%%%%%%%%%%
 \subsection{$\e$-regularization is stable near generic oscillatory interfaces}

As we have seen, the well-posedness of the limit problem at
$\e=0^+$ crucially depends on the asymptotics of the solution
$u(x,t)$ close to the interface point. 

Let us show that, for
typical oscillatory sign-changing behaviour of solutions of the
TFE--4, their $\e$-regularizations do not lead to blow-up.

We restrict our analysis to $N=1$ and assume that, at $t=0$, a
solution profile $u(x,0)$ of \ef{i1} already has a generic form
close to $x=0^+$ being its interface point. Namely, according to
\cite[\S~7]{EGK2}, as $x \to 0^+$,
 \be
 \label{5.13}
  \tex{
 u(x,0) \sim x^\mu [\varphi_*(s)+o(1)] \whereA \mu= \frac 3n,
 \quad s= \ln x,
 }
  \ee
  and the oscillatory component $\varphi_*(s)$ is a periodic solution of the following
  ODE:
 \be
 \label{eqLC}
 \varphi''' + 3(\mu-1) \varphi'' + (3 \mu^2 - 6 \mu +2) \varphi'
+ \mu(\mu-1)(\mu-2)\varphi + |\varphi|^{-n} \varphi=0.
 \ee
 Existence of such a periodic function $\varphi_*(s)$ is known for $n>0$ for some open intervals
(see \cite[\S~7]{EGK2} and more estimates in \cite{PetI})
  \be
  \label{5.14}
   \tex{
  n \in (0,n_{\rm h}) \whereA  n_{\rm h}= 1.7587... < n_+= \frac
  9{3 + \sqrt 3}=1.902...\,,
  }
  \ee
  where $n_+$ is given by the maximal root $\mu_+= \frac 3{n_+}$
  of the quadratic equation (see the third coefficient in the differential operator in \ef{eqLC})
   $$
   \tex{
   3 \mu^2- 6 \mu +2=0 \LongA \mu_{\pm}= \frac {3 \pm \sqrt 3}3.
   }
   $$
   Uniqueness of the periodic $\varphi_*(s)$ remains   an open problem.
  It is worth mentioning that a huge number of numerics  of the ODE \ef{eqLC} always confirmed
  that such a periodic solution is unique for all $n \in (0,n_{\rm h})$ and it is globally stable as $s \to + \iy$,
  i.e., {\em away} from the interface fixed at $y=0^+$.
  It is important that  both of these ``good" properties
  remain valid  up to a heteroclinic bifurcation
   at $n=n_{\rm h}$, at which such a $\varphi_*(s)$ is deformed into a heteroclinic connection of two constant equilibria,
    $$
    \tex{
    \varphi_\pm= \pm \big[\mu_{\rm h}(\mu_{\rm h}-1)(2-\mu_{\rm h})\big]^{\frac 1{n_{\rm h}}}
    \whereA \mu_{\rm h}= \frac 3{n_{\rm h}} \in (1,2)\, .
     }
     $$

    We must admit that
  \ef{5.13} has been obtained from a nonlinear ODE for source-type
  similarity solutions of the TFE--4 in 1D. However, as customary, we expect that
  key evolution properties
  of this ``fundamental" solution remain valid for more general
  solutions of the PDE including the asymptotics \ef{5.13} close
  to the interface. Proving such a generic behaviour for the
  TFE--4 represents a very difficult open problem.

It follows from rescaled data in \ef{5.3} that, for data from  \ef{5.13}, one
requires
 \be
 \label{5.15}
  \tex{
 \a= \frac n3 \andA \b= \frac n3
 }
 \ee
in the group of scaling transformations \ef{5.2}. Overall, this
gives the following rescaled data for $v_\e$:
 \be
 \label{5.16}
  \tex{
 v_{0\e}(y) \sim y^{\frac 3n}[\varphi_*\big( \frac n3 \, \ln \e+ \ln y \big)+ o(1)].
 }
  \ee
Therefore, as $\e \to 0^+$, the data becomes strongly oscillatory in view of fast
phase changes in $\varphi_*(\cdot)$. But what is more important, the data remains uniformly bounded above by a slower,
than in \ef{5.6}, power law:
 \be
 \label{5.17}
 |v_{0\e}(y)| \le {\rm const.} \, |y|^{ \frac 3n} \inB \re.
  \ee
It is not difficult to show that the problem for the rescaled
equation in \ef{5.7} with data as in \ef{5.17} does not admit any
blow-up, and the solution remains bounded on any compact subset.
It is expected that this allows one to match this regular boundary
layer with the outer expansion for $y \gg 1$ to get key features
of solutions $v_\e(y,\t)$ and hence of $u_\e(x,t)$ to see whether
\ef{5.1} actually holds. This is a general strategy towards
uniqueness, and, clearly, much deeper and harder asymptotic
analysis is required here to conclude.

\ssk

Overall, this somehow shows that generic oscillatory behaviour of
source-type is {\em stable} under analytic $\e$-regularizations,
which then cannot lead to any non-uniqueness. However, this is
just a 1D illustration of the regularization procedure. We still
expect that there may exist special complicated configurations of
data in $\ren$ leading to some kind of ``self-focusing"  on
interfaces  implying to non-vanishing $\e$-oscillations as $\e \to
0$ close to interfaces and hence yielding  non-uniqueness. In
general, we expect that, for such higher-order parabolic (and
other) strongly degenerate equations, any global concept of
uniqueness cannot be a key issue in general  PDE theory. Moreover,
for such non-fully-divergent quasilinear degenerate PDEs, an
(absolute) uniqueness can be non-achievable in principle.

%%%%%%%%%%%%%%%%%%%%%%%%%%%%%%%%%%%%%%%%%%%%%%%%%%%%%%%%%%%%%%%%

\end{document}